\documentclass[10pt]{article}

\usepackage{amsmath}
\usepackage{amsfonts}
\usepackage{amssymb}
\usepackage{amsthm}

\newtheorem{theorem}{Theorem}[section]
\newtheorem{lem}{Lemma}[section]
\newtheorem{cor}{Corollary}[section]
\newtheorem{prop}{Proposition}[section]

\newtheorem{klem}{Key Lemma}[section]

\begin{document}

\title{Characterization of the $L^p$-range 
  of the Poisson Transform on the Octonionic Hyperbolic Plane.}
  
\author{ Abdelhamid Boussejra, Nadia Ourchane}

\date{} 
\maketitle
\begin{abstract}
Let $ B(\mathbb{O}^2)=\{x\in \mathbb{O}^2,\mid x\mid<1\}$ be the bounded realization of the exceptional symmetric space $F_{4(-20)}/Spin(9)$. For a nonzero real number $\lambda$, we give a necessary and a sufficient condition on eigenfunctions $F$ of the Laplace-Beltrami operator on $B(\mathbb{O}^2)$ with eigenvalue $-(\lambda^2+\rho^2)$ to have an $L^p$-Poisson integral representations on the boundary $\partial B(\mathbb{O}^2)$. Namely, $F$ is the Poisson integral of an $L^p$-function on the boundary if and only if it satisfies the following growth condition of Hardy-type:
\[
\sup_{0\leq r<1} (1-r^2)^{\frac{-\rho}{2}} \left(\int_{\partial B(\mathbb{O}^2)} |F(r\theta)|^p d\theta\right)^\frac{1}{p}<\infty.
\]
This extends previous results by the first author et al. for classical hyperbolic spaces.\\
\end{abstract}


\section{Introduction}

Let $X=G/K$ be a Riemannian symmetric space of the noncompact type. It is well known that a function $F$ is an eigenfunction of all $G$-invariant differential operators on $X$  if and only if $F$ is  Poisson integral
\[
P_\lambda f(gK)=\int_K f(k)e^{-(i\lambda+\rho)H(g^{-1}k)} dk,
\]
of a hyperfunction $f$ on the Furstenberg boundary $K/M$, for a generic $\lambda\in\mathfrak{a}_\mathbb{C}^*$.\\
( $H$ denotes the projection on the abelian part $A$ of the Iwasawa decomposition of $G$, $\mathfrak{a}_\mathbb{C}^*$ is the complex dual of the Lie algebra $\mathfrak{a}$ of $A$ and $\rho$ is the half sum of the positive roots with multiplicities ).\\
This was conjectured by Helgason  who proved it for the rank one case \cite{H1} and proved in its full generality by M. Kashiwara et al.  \cite{Ka}.\\
A natural question is then  to look for a characterization of the range of the Poisson transform on classical spaces on the Furstenberg boundary $K/M$ such as the spaces $C^{\infty}(K/M), L^p(K/M)$ and the space of distributions $D^\prime(K/M)$. (see \cite{BOS}, \cite{BS}, \cite{Br}, \cite{I}, \cite{Km}, \cite{KS}, \cite{KW}, \cite{Mich}, \cite{L}, \cite{Lo}, \cite{St}, \cite{Sj}). \\
In the case $\lambda=-i\rho$, i.e the harmonic case, the harmonic functions which are Poisson integrals of $L^p$-functions $( 1 < p \leq \infty )$ or bounded measures are characterized by an $H^p$-condition. This was proved by Stoll \cite{St} ( see also Knapp and Williamson \cite{KW} and Michelson \cite{Mich} ).\\
For $\lambda$ in  $\mathfrak{a}^{\ast}_{\mathbb{C}}$ such that $\mathfrak{Re}(i\lambda)$ lies in the open Weyl chamber,  another characterization using weak $L^p$-spaces is given in Sj\"ogren \cite{Sj} ( see also Lohou\'e and Rychener \cite{Lo}).\\
Later on, Ben Said et al. \cite{BOS} gave a characterization of the image $P_{\lambda}(L^p(K/M))$ $( 1<p\leq \infty )$  in terms of Hardy type norm  for $\lambda\in \mathfrak{a}^{\ast}_{\mathbb{C}}$ such that $\mathfrak{Re}(i\lambda)$ lies in the open positive Weyl chamber and the isotropy subgroup of $\lambda$ and $\mathfrak{Im}(\lambda)$ in the Weyl group coincide.\\ 
All the above studies leave out the case $\lambda\in \mathfrak{a}^{\ast}\setminus\{0\}$. Namely the characterization of the image of the Poisson transform on the unitary  spherical principal series representation.\\
Our interest on the problem of characterizing the $L^p$-range of the Poisson transform, for $\lambda\in \mathfrak{a}^{\ast}$  has its root in the work  of Strichartz \cite{S}. More precisely, the statement that, for $\lambda$ in $\mathfrak{a}^{*} \backslash \lbrace 0 \rbrace$ the joint eigenfunctions which are Poisson integrals of $L^2$-functions are characterized by an $L^2$-weight norm, was conjectured by Strichartz in \cite{S} ( for details, see conjecture 4.5 in \cite{S}).\\The Strichartz conjecture in the case  of the complex hyperbolic space, namely $X=SU(n,1)/S(U(n)\times U(1))$, was settled by the first author et al. \cite{BI}. A new proof with an extension to all rank one symmetric spaces was given by Ionescu \cite{I}.  In the case of higher rank the Strichartz conjecture 4.5 was recently solved by K. Kaizuka \cite{K} .\\
In \cite{BS} the first author et al. dealt with Poisson transform of $L^p$-functions. More precisely they proved that Poisson integrals of $L^p$-functions ( $2\leq p <\infty$ ) are characterized by a Hardy type norm in the case $U(n,1,\mathbb{F})/ U(n,\mathbb{F})\times U(1,\mathbb{F})$, $\mathbb{F}=\mathbb{R},\mathbb{C}$ or the quaternions $\mathbb{H}$, i.e the real, complex or quaternionic hyperbolic spaces, respectively (see also Kumar et al.\cite{KS}). The method of the proof uses the techniques of singular integrals on the boundary $K/M$ viewed as a space of homogeneous type in the sense of Coifman and Weiss \cite{C}. Unfortunately, this method depends on the classification of rank one symmetric spaces. That is the classical hyperbolic spaces and the exceptional case.\\
The aim of this paper is to extend  the results in \cite{BI} and \cite{BS} for classical hyperbolic spaces to the case of the Octonionic Hyperbolic plane $F_{4(-20)}/Spin(9)$.\\
In order to describe our result let us fix some notations, referring to Section 2 for more details. \\
Let $\mathbb{O}$ be the division algebra of Octonions ( $\approx$ the Cayley numbers ). Let
\begin{eqnarray*}
B(\mathbb{O}^2)=\{x\in \mathbb{O}^2,\mid x\mid<1\},
\end{eqnarray*}
be the bounded realization of the symmetric space $F_{4(-20)}/Spin(9)$ and let $\partial B(\mathbb{O}^2)$ denote the unit sphere of $\mathbb{O}^2$ with the normalized area measure $d\omega$ on it. Let $L^p(\partial B(\mathbb{O}^2))$ denote  the space of all $\mathbb{C}$-valued measurable (classes) functions $f$ on $\partial B(\mathbb{O}^2)$ with $\Vert f \Vert_p<\infty$. Here  
\begin{eqnarray*}
\Vert f \Vert_p=\left( \int_{\partial B(\mathbb{O}^2)} |f(\omega)|^p d\omega \right)^{\frac{1}{p}}.
\end{eqnarray*}
For $f\in L^p(\partial B(\mathbb{O}^2))$ and $\lambda$ a complex number, we define  the Poisson transform of $f$ by 
\begin{eqnarray*}
P_\lambda f(x)=\int_{\partial B(\mathbb{O}^2)}\left(\frac{1-\mid x\mid^2}{\mid 1-[x,\omega]\mid^{2}}\right)^{\frac{i\lambda+\rho}{2}}f(\omega)d\omega,\quad \rho=11.
\end{eqnarray*}
For the precise definition of $[x,\omega]$ see (\ref{[x,y]}) section 2.\\
Let $\mathcal{E}_{\lambda}(B(\mathbb{O}^2))$ be the space of all eigenfunctions of  the Laplace-Beltrami operator  $\Delta$ of $B(\mathbb{O}^2)$ with eigenvalue $-(\lambda^2+\rho^{2})$.\\
In order to characterize, for $\lambda\in\mathbb{R}\setminus\{0\}$, those $F\in \mathcal{E}_{\lambda}(B(\mathbb{O}^2))$ which are Poisson transform by $P_\lambda$ of some $f$ in $L^p(\partial B(\mathbb{O}^2))$ ( $1< p < \infty $ ), we introduce the Hardy type space $\mathcal{E}_{\lambda,p}^*(B(\mathbb{O}^2))$  consisting of functions $F\in \mathcal{E}_{\lambda}(B(\mathbb{O}^2))$ such that
\begin{eqnarray*}
 \|F\|_{*,p}=\sup_{0\leq r<1} (1-r^2)^{\frac{-\rho}{2}} \left(\int_{\partial B(\mathbb{O}^2)} |F(r\theta)|^p d\theta\right)^\frac{1}{p}<\infty.
\end{eqnarray*} 
Finally, we denote by $B(o,t)$ the geodesic ball of radius $t$ centred at $0$ and $d\mu(x)$ the $G$-invariant measure of $B(\mathbb{O}^2)$, given by $d\mu(x)=(1-\mid x\mid^2)^{-\rho-1}dm(x)$, $dm(x)$ being the Lebesgue measure.\\In the sequel of this paper  $c $ will denote a numerical positive constant and $\rho=11$.\\
The first main result we prove in this paper can be stated as follows:
\begin{theorem}\label{th1}
Let $\lambda \in \mathbb{R}\setminus \{0\} $. Then we have
\begin{itemize}
\item[i)] A function $F\in\mathcal{E}_{\lambda}(B(\mathbb{O}^2)) $ is the Poisson transform by $P_{\lambda}$ of some $f$ in $L^2(\partial B(\mathbb{O}^2))$  if and only if $F\in \mathcal{E}_{\lambda,2}^*(B(\mathbb{O}^2))$.\\
Moreover, there exists a positive constant  $c $ such that for every $f\in L^2(\partial B(\mathbb{O}^2))$ the following estimates hold:
\begin{equation}\label{estimate1}
|c(\lambda)|\|f\|_2\leq \|P_{\lambda}f\|_{*,2} \leq c ( 1+\mid\lambda\mid+\frac{1}{\mid\lambda\mid} )\|f\|_2.
\end{equation}
\item[(ii)] Let $F\in \mathcal{E}_{\lambda,2}^*(B(\mathbb{O}^2))$. Then its $L^2$-boundary value $f$ is given by \\the following inversion formula:
\begin{equation}\label{invers}
f(\omega)=\mid c(\lambda)\mid^{-2}\lim_{t\rightarrow +\infty}\frac{1}{t}\int_{B(o,t)}P_{-\lambda}(x,\omega)F(x)d\mu(x),
\end{equation} 
in $L^2(\partial B(\mathbb{O}^2))$.
\end{itemize}
\end{theorem}
In (\ref{estimate1}), $c(\lambda)$ is the Harish-Chandra c-function associated to the Octonionic hyperbolic plane $B(\mathbb{O}^2)$, given by 
\begin{eqnarray*}
c(\lambda)=\frac{\Gamma(8)\Gamma(i\lambda)}{\Gamma(\frac{i\lambda+\rho}{2}-3)\Gamma(\frac{i\lambda+\rho}{2})}.
\end{eqnarray*}
The second major result of this paper is
\begin{theorem}\label{th2}
Let $\lambda \in \mathbb{R}\setminus \{0\} $ and let $p\in ]1,+\infty[$. Then we have, a function $F\in\mathcal{E}_{\lambda}(B(\mathbb{O}^2)) $ is the Poisson transform by $P_{\lambda}$ of some $f\in L^p(\partial B(\mathbb{O}^2))$  if and only if $F\in \mathcal{E}_{\lambda,p}^*(B(\mathbb{O}^2))$.\\
Moreover, there exists a positive constant  $\gamma(\lambda,p)$ such that for every $f\in L^p(\partial B(\mathbb{O}^2))$ the following estimates hold:
\begin{equation}\label{estimate2}
|c(\lambda)|\|f\|_p\leq \|P_{\lambda}f\|_{*,p} \leq \gamma(\lambda,p)\|f\|_p.
\end{equation}
\end{theorem}
As mentioned in the introduction, the main difficulty in proving our results lies in proving some uniform pointwise estimates on the generalized spherical functions which unlike the case $\lambda\in \mathbb{C}\setminus\mathbb{R}$ have oscillating terms at infinity.\\
In this paper we prove  the right hand side of the estimate (\ref{estimate1}) in  Theorem \ref{th1} by adapting to the Octonionic case the method that we used  in the case of the classical hyperbolic spaces, see \cite{BS}. More precisely, we will discuss a uniform $L^p$-boundedness of a family of Calderon-Zygmund operators $(\Psi_r(\lambda))_{r\in [0,1[}$ ( see Section 4) on the boundary  $\partial B(\mathbb{O}^2)$ considered as a space of homogeneous type in the sense of Coifman and Weiss \cite{C}. 
To prove the sufficiency condition for $p\neq 2$, we follow the method we used in \cite{B1} and \cite{B2}, to characterize Poisson integrals of $L^p$-functions on the Shilov boundary of bounded symmetric domains.\\Although  the techniques we use here may seem to be similar to those in \cite{BS}, however in working in the exceptional case we encounter a prime difficulty, due to the fact that the algebra of Octonions $\mathbb{O}$ is not associative.
\\
\textbf{Consequences.}\\
(i) As a direct consequence of Theorem \ref{th2}, we obtain that for $\lambda\in \mathbb{R}\setminus \{0\}$ and $p\in ]1,\infty[$, the Hardy-type spaces $\mathcal{E}_{\lambda,p}^*(B(\mathbb{O}^2))$ are Banach spaces. This is closely similar to O.Bray conjecture in the case of the real hyperbolic space with $p=2$, see \cite{Br}.\\
(ii) Let $\Phi_{\lambda,lm}$ be the generalized spherical function associated to the Octonionic hyperbolic plane (see Proposition \ref{prop1}). Namely
\begin{eqnarray*}
\begin{split}
\Phi_{\lambda,lm}(r)&=(8)_{l}^{-1}(\frac{i\lambda+\rho}{2})_{\frac{m+l}{2}}(\frac{i\lambda+\rho}{2}-3)_{\frac{l-m}{2}}r^{l}(1-r^{2})^{\frac{i\lambda+\rho}{2}}\\
&_2F_1(\frac{i\lambda+\rho+l+m}{2},\frac{i\lambda+\rho+l-m}{2}-3,l+8;r^2),
\end{split}
\end{eqnarray*}
with $(a)_k=a(a+1)...(a+k-1)$ is the Pochammer symbol and  $_2F_1(a,b,c;x)$ is the Gauss hypergeometric function.\\
Then as an immediate consequence of Theorem \ref{th1} we obtain the following uniform pointwise estimate. 
\begin{cor}
Let $\lambda$ be a nonzero real number. Then there exists a positive constant $c $ such that 
\begin{eqnarray*}
\sup_{(l,m)\in \widehat{K}_{0}}\mid \Phi_{\lambda,lm}(r)\mid\le c (1+\mid\lambda\mid+\frac{1}{\mid\lambda\mid})((1-r^2)^{\frac{\rho}{2}}, 
\end{eqnarray*}
for all $r\in [0,1[$. 
\end{cor}
As $\Phi_{\lambda,lm}$ can be written in terms of Jacobi functions, the above estimate might have independent interest on its own from a view point of special functions.\\
(iii) Recently in their study of Roe Theorem in rank one symmetric spaces,\\ Kumar et al.\cite{KS}  proved a characterization of the $L^p$-range of the Poisson transform by means of Lorentz type norms in the case of the classical hyperbolic spaces $B(\mathbb{F}^n)$. Below we will show how to apply  Theorem \ref{th2}  to extend their result to the exceptional case.\\
We follow mainly the notations in \cite{KS}. Let $X=G/K$ be a hyperbolic space on $\mathbb{R}, \mathbb{C}$ or $\mathbb{H}$. Let $\mu$ denote the $G$-invariant measure of $X$. For $F$ a $\mu$-measurable  complex-valued function on $X$, we set $\parallel F\parallel_{p,\infty}=\sup_{s>0}sd_{F}(s)^\frac{1}{p}$, where $d_{F}$ is the distribution function of $F$. Finally  define $\mathcal{A}_{2,p}(F)=\Vert \mathcal{A}_p(F) \Vert_{2,\infty}$, with $\mathcal{A}_p(F)(x)=\left( \int_{K/M}|F(kx)|^p dk \right)^{\frac{1}{p}}$ ( where as usual $M$ is the centralize of $A$ in $K$, if $G=KAN$ is an Iwasawa decomposition of $G$ ).\\
Then Kumar et al. result may be state as follows
\begin{theorem}[\cite{KS}]
Let $1<p<\infty$ and $X$ a hyperbolic space over $\mathbb{R}, \mathbb{C}$ or $\mathbb{H}$. Let $F$ be a $\mathbb{C}$-valued function on $X$. Then $F$ is the Poisson transform by $P_{\lambda}$ of some $f\in L^p(\partial B(\mathbb{F}^n))$  if and only if it satisfies $\mathcal{A}_{2,p}(F)<\infty$.
\end{theorem}
As a consequence of the method of the proof of Theorem \ref{th2}, we extend the above result to the case $X$ is the Octonionic hyperbolic space.
\begin{theorem}
Let $\lambda \in \mathbb{R}\setminus \{0\} $ and $p\in ]1,+\infty[$. Let $F$ be  a $\mathbb{C}$-valued function on the Octonionic hyperbolic plane $B(\mathbb{O}^2)$, satisfying $\Delta F=-(\lambda^2+\rho^2)F$. Then $F=P_{\lambda}f$, for some $f\in L^{p}(\partial B(\mathbb{O}^2))$ if and only if $\mathcal{A}_{2,p}(F)<\infty$.
\end{theorem}
\textbf{Proof}. Noting that if $\|F\|_{\ast,p}<\infty$, then $\mathcal{A}_{2,p}(F)\leq c \|F\|_{\ast,p} $, the necessary condition follows from the right hand side of the estimate (\ref{estimate2}) in Theorem \ref{th2}. Next, using the $L^2$ inversion formula (\ref{invers}), we may follow the same method as in the proof of the sufficiency condition of Theorem 4.3.6 in \cite{KS} to get the result. So we omit it.\\
Now we give the organization of this paper. In section 2, some preliminaries of harmonic analysis on the Cayley plane are described. In section 3 we prove our main results. The proof rely on  establishing  the Key Lemma of this paper giving a uniform $L^2$-boundedness of the Calderon-Zygmund operators $(\Psi_{r}(\lambda))$ ($r\in[0,1[$) associated to the Poisson transform $P_{\lambda}$. To prove the Key Lemma, we will adapt on $\partial B(\mathbb{O}^2)$ - in a uniform manner in $r\in[0,1[$ - the method of proving the $T(1)-$ Theorem of David-Journ\'e and Semens \cite{D} and for this we follow the program accomplished by Y.Meyer in his new proof of the $T(1)$-Theorem in $L^{2}(\mathbb{R}^{n})$ which is based on the Cotlar-Stein Lemma. This is the subject of section 4.\\

We end this section with a brief discussion on all rank one symmetric spaces. Let $X=G/K$ be a noncompact Riemannian symmetric space of rank one. Then $X$ can be  realized  as  the unit ball in $\mathbb{F}^n$, with $\mathbb{F}$ is either the real numbers $\mathbb{R}$, or the complex numbers $\mathbb{C}$, or the quaternionic $\mathbb{H}$, or the Cayley numbers $\mathbb{O}$, in last case $n=2$. Moreover,\\ 
if $\mathbb{F}= \mathbb{R}$, then $(G,K)=(SO_e(n,1),SO(n))$\\
if $\mathbb{F}= \mathbb{C}$, then $(G,K)=(SU(n,1),S(U(n)\times U(1)))$\\
if $\mathbb{F}= \mathbb{H}$, then $(G,K)=(Sp(n,1),Sp(n)\times Sp(1))$\\
if $\mathbb{F}= \mathbb{O}$, then $n=2$ and $(G,K)=(F_{4(-20)},Spin(9))$.\\
Let $\rho_{\mathbb{F}}$ denote the half sum of positive roots of the pair $(\mathfrak{g},\mathfrak{a})$, $\mathfrak{g}$ being the Lie algebra of $G$ and $\mathfrak{a}$ a Cartan subalgebra of $\mathfrak{g}$ ( $\mathfrak{a}=\mathbb{R}H_0$, since $rank(X)=1$ ). By abuse of notation we will denote $\rho_{\mathbb{F}}(H_0)$ by $\rho_{\mathbb{F}}$. Then $\rho_{\mathbb{F}}=\frac{n-1}{2}, n, 2n+1$ or $11$ accordingly to $\mathbb{F}=\mathbb{R}, \mathbb{C}, \mathbb{H}$ or $\mathbb{O}$.\\ 
Denote by $A$ the analytic subgroup of $G$ that corresponds to $\mathfrak{a}$. Then $A$ may be parametrized by $a_t=\exp tH_0$.\\
For $p\in]1,\infty[$, let $M^{\ast}_p(F)=\sup_{t>0}e^{\rho t}\left(\int_{K}\mid F(ka_t.0)\mid^p dk\right)^{\frac{1}{p}}$, where $dk$ is the normalized Haar measure of $K$.\\
Now, taking into account  Theorem \ref{th1}, Theorem \ref{th2} and the results in  ( \cite{BS}, Theorem A ) and (\cite{KS}, Theorem 4.3.6), a characterization of the $L^p$-range of the Poisson transform  for the rank one symmetric spaces is now completed. With the help of the above notations theses results may be state in a unified manner as follows:
\begin{theorem}\label{th3}
Let $\lambda\in \mathfrak{a}^\ast\setminus\{0\}$ and let $1<p<\infty$. Let $F$ be a $\mathbb{C}$-valued function on $X$ satisfying $\Delta F=-(\lambda^2+\rho^2)F$. Then we have
\begin{itemize}
\item[i)]$F$ has an $L^p$-Poisson integral representation on $K/M$ if and only if \\ $M^{\ast}_p(F)<\infty$. 
Moreover, there exists a positive constant  $\gamma(\lambda,p)$ such that for every $f\in L^p(K/M)$ the following estimates hold:
\[
|c(\lambda)|\|f\|_{p}\leq M^{\ast}_p(P_\lambda f) \leq \gamma(\lambda,p)\|f\|_{p}.
\]
\item[ii)]If $M^{\ast}_2(F)<\infty$. Then $F$ has an $L^2$-boundary value $f$ given by the following inversion formula
\[
f(kM)=\lim_{t\rightarrow +\infty}\frac{1}{t}\int_{B(0,t)}F(x)P_{-\lambda}(x,kM)d\mu(x)
\]
in $L^{2}(K/M)$. 
\end{itemize}
\end{theorem}
In above  $B(0,t)$ denotes the geodesic ball of radius $t$ centred at $0$ and the $G$-invariant measure $d\mu$ on $X$.\\  
Below we give the last result of this section. \\
For any locally integrable function $F$ with respect to $d\mu$, we set  
\[
M_2(F)^{2}=\sup_{t>0}\frac{1}{t}\int_{B(0,t)}\mid F(x)\mid^2 d\mu(x).
\]
Writing $x=ka_t.0$, we easily see that if $M^{\ast}_2(F)$ is finite then 
\begin{equation}\label{H}
M_2(F)\leq c M^{\ast}_2(F),
\end{equation}
for some $c $ positive constant. Therefore if  $f\in L^2(\partial B(\mathbb{F}^n))$, then 
\[
M_2(P_\lambda f)\leq c (1+\mid \lambda\mid+\frac{1}{\mid\lambda\mid})\parallel f\parallel_{2},
\]
by the above Theorem. Thus  we have proved 
\begin{cor}
Let $\lambda \in \mathbb{R}\setminus \{0\} $ and let $f\in L^2(\partial B(\mathbb{F}^n))$. Then
\begin{eqnarray*}
\left(\sup_{t>0}\frac{1}{t}\int_{B(0,t)}\mid P_{\lambda}f(x)\mid^2 d\mu(x)\right)^{\frac{1}{2}}\le c (1+\mid \lambda\mid+\frac{1}{\mid\lambda\mid})\parallel f\parallel_{2}.
\end{eqnarray*}
\end{cor}
The above result  has already been established by Ionescu using different method, see \cite{I}.

\section{Preliminary results}
We review in this section some known results of harmonic analysis on the Octonionic hyperbolic space. We first recall some  properties on the Octonions that will be needed in this paper, referring to  \cite{T} and \cite{V} for more details.\\
We denote by $\mathbb{O}$ the algebra of Octonions. $\mathbb{O}$ has a basis over $\mathbb{R}$ given by  $e_0, e_1,.....,e_7$, where $e_0$ is the unit element, and $e_m$ are anti-commuting elements satisfying $e_m^2=-1$.\\ 
We define the standard involution of $\mathbb{O}$ over $\mathbb{R}$ by $\overline{x}=x_0-\sum\limits^7_{j=1} x_je_j$, and we have $\overline{xy}=(\overline{y})(\overline{x})$, for every $x,y\in \mathbb{O}$. \\
If $x=\sum\limits_{j=0}^{7} x_je_j$, the summand  $x_0 e_0=x_0$ is called the real part of $x$ and it is noted by $\mathfrak{Re}(x)$. Furthermore, the norm  on $\mathbb{O}$ is defined as $\mid x\mid^2=\sum\limits^7_{j=0} x^{2}_j$, and it satisfies the identity $\mid xy\mid=\mid x\mid \mid y\mid$. Every nonzero octonion $x$ has a unique inverse, namely
\begin{eqnarray*}
x^{-1}=\mid x\mid^{-2}\overline{x}.
\end{eqnarray*}
For $x,y\in \mathbb{O}^2$, we put
\[
\Phi(x,y)=\sum_{j=1}^{2}\mid x_j\mid^2\mid y_j\mid^2+2\Re((x_1x_2)(\overline{y_1y_2})).
\]
The form $\Phi(x,y)$ may be written as 
\begin{equation}\label{Phi(x,y)=||}
\Phi(x,y)=\mid (\overline{x_1}y_2)(y_2^{-1}y_1)+x_2\overline{y_2}\mid^{2},
\end{equation}
for $y_2\neq 0$.\\
Also,  we consider $O_{\mathbb{O}}(2)$ the group of all $\mathbb{R}-$linear transformations of $\mathbb{R}^{16}$ which preserve the form $\Phi(x,y)$. The group $O_{\mathbb{O}}(2)$ is  a subgroup of the  orthogonal group $O(16)$  ( when identifying $\mathbb{O}^2$ with $\mathbb{R}^{16}$ ).\\
Let $\partial B(\mathbb{O}^2)=\{\omega\in \mathbb{O}^2: \mid \omega\mid=1\}$ be the unit sphere in $\mathbb{O}^2$. Then we have 
\begin{lem}[\cite{V}]
The group $O_{\mathbb{O}}(2)$ acts transitively on the unit sphere $\partial B(\mathbb{O}^2)$.
\end{lem} 
\subsection{\textbf{Bounded realization of $F_{4(-20)}/Spin(9)$.}}
In this section we give, after Takahashi\cite{T}, treatments of  the Octonionic hyperbolic space that  are  suitable for computation and that are helpful to handle the uniform $L^p-$characterization of the Poisson transform on $F_{4(-20)}/Spin(9)$ \\( see also \cite{V} ).\\  
Let $A(3,\mathbb{O})$ be the exceptional Jordan algebra of $3 \times 3 $ Hermitian matrices with entries from $\mathbb{O}$ and let $A(3,\mathbb{O} \otimes \mathbb{C})$ be its complexification. The Jordan  product  being $A\circ B=\frac{1}{2}(AB+BA)$.\\
Let $J$ denote the Jordan subalgebra of $A(3,\mathbb{O} \otimes \mathbb{C})$, consisting of matrices 
 \begin{eqnarray*}
 X = \left(
\begin{array}{ccc}
a_{1} & u_{3}\otimes (-1)^{\frac{1}{2}} & \overline{u_{2}}\otimes (-1)^{\frac{1}{2}} \\
\overline{u_{3}}\otimes (-1)^{\frac{1}{2}} & a_{2} & u_{1} \\
u_{2}\otimes (-1)^{\frac{1}{2}} & \overline{u_{1}} & a_{3}
\end{array}
\right),
 \end{eqnarray*}
with $a_j \in \mathbb{R}$ and $u_j \in \mathbb{O}$, for $j=1,2,3$.\\
We denote by $\Omega$ the subset of $J$ consisting of idempotents elements  $X\in J$ such that $tr(X)=1$ and $tr(X\circ E_1)\ge 1$. Here  
 $E_1=\left(
\begin{array}{ccc}
1 & 0 & 0 \\
0 & 0 & 0 \\
0 & 0 & 0
\end{array}
\right)$.\\
Next, let $B(\mathbb{O}^2)=\lbrace (x_1,x_2) \in \mathbb{O}^2, |x_1|^2+|x_2|^2<1  \rbrace$ be the unit ball of $\mathbb{O}^2$. Then the map 
\begin{eqnarray*}
(x_1,x_2) \longmapsto X(x_1,x_2)=\dfrac{1}{1-|x|^2}\left(
\begin{array}{ccc}
1 & \overline{x_2}\otimes (-1)^{\frac{1}{2}} & \overline{x_1}\otimes (-1)^{\frac{1}{2}} \\
x_2\otimes (-1)^{\frac{1}{2}} & -|x_2|^2 & -x_2\overline{x_1} \\
x_1\otimes (-1)^{\frac{1}{2}} & -x_1\overline{x_2} & -|x_1|^2
\end{array}
\right)
\end{eqnarray*}
is a bijection from the unit ball  $B(\mathbb{O}^2)$ onto $\Omega$.\\
Let $G=F_{4(-20)}$ be the identity component of $Aut(J)$ the group of all automorphisms of $J$. Then using the above map we shall define an action of $G$ on $B(\mathbb{O}^2)$.\\
Namely for $x=(x_1,x_2)\in B(\mathbb{O}^2)$, we set 
$g \cdot (x_1,x_2)=(x'_1,x'_2)$ according to $gX(x_1,x_2)=X(x'_1,x'_2)$.\\
The action of $G$ on $B(\mathbb{O}^2)$ is transitive and as  homogeneous space $G \slash K$ is then identified to $B(\mathbb{O}^2)$, where $K$ is the isotropy subgroup of $E_1$. Moreover $K$ is a maximal compact subgroup of $G$ which is isomorphic to $Spin(9)$.\\
Next, we identify the boundary $\partial B(\mathbb{O}^2)$ to the set 
\begin{eqnarray*}
\lbrace Y(u,v)=\left(
\begin{array}{ccc}
0 & v & \overline{u} \\
\overline{v} & 0 & 0 \\
u & 0 & 0
\end{array}
\right): u,v \in \mathbb{O}, |u|^2+|v|^2=1\rbrace.
\end{eqnarray*}
Then using this identification, the action of $K$  on $\partial B(\mathbb{O}^2)$ is defined as follows:\\ 
For $u,v$ in $\partial B(\mathbb{O}^2)$ put $k.(u,v)=(u^{\prime},v^{\prime})$, via $kY(u,v)=Y(u^{\prime},v^{\prime})$.
\\ 
It is shown in \cite{T} that this action is transitive  and as homogeneous space we have $\partial B(\mathbb{O}^2) = K\slash M$, where $M$ is the isotropic subgroup of the element $F_2^1=\left(
\begin{array}{ccc}
0 & 0 & 1 \\
0 & 0 & 0 \\
1 & 0 & 0
\end{array}
\right)$. Moreover the group $M$ is isomorphic to $Spin(7)$.\\
Now we describe, in brief, the Peter-Weyl decomposition of $L^2(\partial B(\mathbb{O}^2))$ under the action of $K$.
Let $\widehat{K}_{0}$ be the set of  pairs $(l,m)$ of nonnegative integers  such that $l\geq m\geq 0$ and $l\pm m$ is even. Then under the action of $K$, the Peter-Weyl decomposition of $L^2(\partial B(\mathbb{O}^2))$ is given by 
\begin{eqnarray*}
L^2(\partial B(\mathbb{O}^2))=\bigoplus_{(l,m)\in \widehat{K}_{0}}V^{lm},
\end{eqnarray*}
where $V^{lm}$ is the $K-$cyclic space for the zonal spherical function associated to the pair $(K,M)$.(See \cite{T} and \cite{J} ).
\subsection{\textbf{The Poisson transform.}}
In this subsection we recall some known results on the Poisson transform $P_\lambda$ associated to the Octonionic hyperbolic plane $B(\mathbb{O}^2)$ that will be needed in the sequel.\\
Firstly, we introduce some notations. For $x,y\in \mathbb{O}^{2}$, put 
\begin{equation}\label{Psi(x,y)}
\Psi(x,y)=1-2<x,y>_{\mathbb{R}}+\Phi(x,y),
\end{equation}
where $<x,y>_{\mathbb{R}}$ is the Euclidean inner product of $x$ and $y$ as vectors in $\mathbb{R}^{16}$.\\
Also, for $x,y\in \mathbb{O}^2$, we set 
\begin{equation}\label{[x,y]}
[x,y]=
\left\{
\begin{array}{rl}
(\overline{x_1}y_2)(y_2^{-1}y_1)+x_2\overline{y_2}& \ if \ y_2\neq0 \\
\overline{x_1}y_1 & \ if \ y_2 =0.
\end{array}
\right.
\end{equation}
Then $\Psi$ may be written as 
\begin{equation}\label{Psi=1-[]}
\Psi(x,y)=\mid 1-[x,y]\mid^{2},
\end{equation}
and since $|[x,y]|\leq |x||y|$, we easily see that $\Psi(x,y)>0$, for $x\in\overline{B(\mathbb{O}^2)}$ and $y\in  B(\mathbb{O}^2)$.\\
Note that $\mid[x,y]\mid$ is the analogue of the form $\mid \sum\limits_{j=1}^{n}a_j\overline{b_j}\mid$ in the real, complex or the quaternionic fields.\\  
The Poisson kernel for the Octonionic hyperbolic plane is the function $P(x,\omega)$ defined on $B(\mathbb{O}^2)\times\partial B(\mathbb{O}^2)$ by
\begin{eqnarray*}
P(x,\omega)=\left(\frac{1-\mid x\mid^2}{\Psi(x,\omega)}\right)^{\rho},\quad \rho=11
\end{eqnarray*}
For $f\in L^1(\partial B(\mathbb{O}^2))$ and $\lambda\in \mathbb{C}$, we define the Poisson integral of $f$ by 
\begin{eqnarray*}
P_{\lambda}f(x)=\int_{\partial B(\mathbb{O}^2)}P_{\lambda}(x,\omega)f(\omega)d\omega,
\end{eqnarray*} 
where $P_{\lambda}(x,\omega)=\left[ P(x,\omega)\right]^{\frac{i\lambda+\rho}{2\rho}}$.\\
Below we recall a result on the precise action of the Poisson transform on the $K-$types $V^{lm}$. 
\begin{prop}[\cite{H1}]\label{prop1}
Let $\lambda$ be  a complex number and let $f\in V^{lm}$. Then we have 
\begin{eqnarray*}
P_{\lambda}f(r\theta)=\Phi_{\lambda,lm}(r)f(\theta),
\end{eqnarray*}
where $\Phi_{\lambda,lm}(r)$ is the generalized spherical function given by
\begin{eqnarray*}
\begin{split}
\Phi_{\lambda,lm}(r)&=(8)_{l}^{-1}(\frac{i\lambda+\rho}{2})_{\frac{m+l}{2}}(\frac{i\lambda+\rho}{2}-3)_{\frac{l-m}{2}}r^{l}(1-r^{2})^{\frac{i\lambda+\rho}{2}}\\
&_2F_1(\frac{i\lambda+\rho+l+m}{2},\frac{i\lambda+\rho+l-m}{2}-3,l+8;r^2).
\end{split}
\end{eqnarray*}
\end{prop}
We end this section by a result on the asymptotic behaviour of the generalized spherical functions, which is due to Ionescu\cite{I}.
\begin{lem}\label{lem2}
Let $\lambda$ be a nonzero real number. Then there exists a positive constant $c $ such that 
\begin{equation}
\lim_{t\rightarrow +\infty}\frac{1}{t}\int_{B(0,t)}\mid \Phi_{\lambda,l m}(\mid x\mid)\mid^{2}d\mu(x)=c \mid c(\lambda)\mid^2,
\end{equation}
for every $(l,m) \in \widehat{K}_{0}$.
\end{lem}
\section{Proof of the Main Results}
As explained in the introduction, the main difficulty in proving Theorem \ref{th1} and Theorem\ref{th2} is to show that the image $P_\lambda(L^p(\partial B(\mathbb{O}^2)$ is continuously embedded in $\mathcal{E}_{\lambda,p}^*(B(\mathbb{O}^2))$. We will overcome  this difficulty  by discussing the uniform $L^{p}-$boundedness of the following  family of superficial Poisson-Szeg\"o integrals $(\Psi_{r}(\lambda))_{r\in [0,1[}$ :
\begin{equation}
\Psi_{r}(\lambda)f(\theta)=\int_{\partial B(\mathbb{O}^2)}\Psi_{r}(\lambda,\theta,\omega)f(\omega)d\omega,
\end{equation}
where the Schwartz kernel is given by $\Psi_{r}(\lambda,\theta,\omega)=\left(\Psi(r\theta,\omega)\right)^{\frac{-i\lambda-\rho}{2}}$.\\
To do so, we equip the unit sphere $\partial B(\mathbb{O}^2)$ with the following non-isotropic metric $d(a,b)=\mid 1-[a,b]\mid^{\frac{1}{2}}$ ( see Section 4 ), so that $(\partial B(\mathbb{O}^2),d)$ becomes a space of homogeneous type in the sense of Coifman and Weiss. Now we state the Key Lemma of this paper.
\begin{klem}\label{klem}
Let $\lambda$ be a nonzero real number. Then there exists a positive constant $c $ such that the following estimates hold
\begin{equation}\label{Q1}
\quad \sup_{0\leq r<1} ||\Psi_{r}(\lambda)||_2 \leq c (1+\mid \lambda\mid+\frac{1}{\mid \lambda\mid}),
\end{equation}
where $||.||_2$ stands for the $L^2-$operatorial norm.
\begin{equation}\label{Q2}
\quad \sup_{0\leq r<1}\int_{d(\omega,e_1)>2d(\theta,e_1)}\mid \Psi_{r}(\lambda,\omega,\theta)-\Psi_{r}(\lambda,\omega,e_1)\mid d\omega\le c (1+\mid \lambda\mid).
\end{equation}
\end{klem}
Notice that from the estimates (\ref{Q1}) as well as the H\"ormander condition  (\ref{Q2}) we deduce that the following estimate holds  
\begin{equation}\label{Q3}
\sup_{0\leq r<1}|| \Psi_{r}(\lambda)||_p \leq \gamma(\lambda,p),\quad \textit{for} \ p\in ]1,\infty[,
\end{equation}
by the Marcinkiewicz interpolation theorem and duality.\\
\textbf{Proof of Theorem \ref{th1}.} \\
i)The necessary condition. Let $f\in L^2(\partial B(\mathbb{O}^2)$, since
\begin{equation}\label{Poisson}
P_{\lambda}f(r\theta)=(1-r^2)^{\frac{i\lambda+\rho}{2}}\Psi_{r}(\lambda)f(\theta),
\end{equation}
then  
\begin{eqnarray*}
\parallel P_{\lambda}f\parallel_{\ast,2}\le c(1+\mid\lambda\mid+\frac{1}{\mid\lambda\mid}\|f\|_2,
\end{eqnarray*} 
by (\ref{Q1}). This prove the right hand side of the estimate (\ref{estimate1}) in Theorem \ref{th1}.\\
To prove the sufficiency condition, let $F\in \mathcal{E}_{\lambda}(B(\mathbb{O}^2))$ satisfying the growth condition $\parallel F\parallel_{\ast,2}<\infty$. Recall that the $G-$invariant measure of $B(\mathbb{O}^2)$ is given by $d\mu(x)=(1-\mid x\mid^2)^{-\rho-1}dm(x)$, $dm(x)$ being the Lebesgue measure. \\Let $M_2(F)$ denote the following $L^2-$weighted norm 
\[
M_2(F)^2=\frac{1}{t}\int_{B(0,t)}\mid F(x)\mid^2d\mu(x).
\]
Using the polar coordinates $x=r\theta$, we easily see that if $\parallel F\parallel_{\ast,2}$ is finite, then 
\[
M_2(F)\le c  \parallel F\parallel_{\ast,2},
\]
for some positive constant $c$.\\
Next,  by  Ionescu result \cite{I}, we know that there exists $f\in L^{2}(\partial B(\mathbb{O}^2))$ such that $F=P_{\lambda}f$ with 
\[
\mid c(\lambda)\mid \parallel f\parallel_{2} \le M_2(P_{\lambda}f).
\] 
Therefore $ \mid c(\lambda)\mid \parallel f\parallel_{2} \le \parallel F\parallel_{\ast,2}$ and the proof of i) is finished.\\
ii)  Although the proof of the inversion formula is the same as in the classical hyperbolic spaces \cite{BS}, but for the sake of completeness we give here the outline of the proof.\\ 
Let $F\in \mathcal{E}_{\lambda}(B(\mathbb{O}^2))$ such that $\parallel F\parallel_{*,2}<\infty$. Then $F=P_{\lambda}f$, for some $f$ in  $L^{2}(\partial B(\mathbb{O}^2))$. Expanding $f$ into its $K-$type series, $f=\sum_{(l,m)\in \widehat{K}_{0}}f_{l m}$  and using Proposition \ref{prop1}, $F$ may be written as   
\begin{eqnarray*}
F(r\theta)=\sum\limits_{(l,m)\in \widehat{K}_{0}}\Phi_{\lambda,l m}(r)f_{lm}(\theta),
\end{eqnarray*}
in $C^{\infty}([0,1[\times\partial\mathbb{B}(\mathbb{O}^2))$.\\
Next, set 
\begin{eqnarray*}
g_t(\omega)=\mid c(\lambda)\mid^{-2}\frac{1}{t}\int_{B(0,t)}F(x)P_{-\lambda}(x,\omega)dx.
\end{eqnarray*}
Then replacing $F$ by  its  expansion series and using Proposition \ref{prop1} once again, we get 
\begin{eqnarray*}
g_t(\omega)=\frac{1}{t}\sum_{(l,m)\in \widehat{K}_{0}}\left( \int_{0}^{\tanh(t)}\mid \Phi_{\lambda,lm}(r)\mid^{2}(1-r^2)^{-\rho-1} r^{15}dr \right) f_{l m}(\omega).
\end{eqnarray*}
It is easy to see that the functions $(g_t)_{t>0}$ are in $L^{2}(\partial B(\mathbb{O}^2))$. Furthermore 
\begin{eqnarray*}
\parallel g_t-f\parallel_2^{2}=\sum_{(l,m)\in \widehat{K}_{0}}\left| \frac{\mid c(\lambda\mid^{-2}}{t}\int_0^{\tanh(t)}|\Phi_{\lambda,lm}(r)|^{2}(1-r^2)^{-\rho-1} r^{15}dr - 1 \right|^2\parallel f_{lm}\parallel_2^2.
\end{eqnarray*}
Next, using the uniform asymptotic behaviour of the generalized spherical functions $\Phi_{\lambda,lm}(r)$ in Lemma \ref{lem2}, we get 
\begin{eqnarray*}
\lim_{t\rightarrow \infty}\parallel g_t-f\parallel_2^{2}=0.
\end{eqnarray*}
This finishes the proof of Theorem \ref{th1}.\\
\textbf{Proof of Theorem \ref{th2}}\\
Let $f\in L^p(\partial B(\mathbb{O}^2)$. Then using the identity (\ref{Poisson}) as well as the estimate (\ref{Q3}), we get 
\[
\parallel P_{\lambda}f\parallel_{\ast,p}\le \gamma(\lambda,p)\|f\|_p.
\]
This prove the necessary condition in Theorem \ref{th2}.\\
Next, to prove the sufficiency condition of Theorem \ref{th2} let  $(\chi_n)_n$ denote an approximation of the  identity in $ \mathcal{C}(K)$. That is 
\begin{eqnarray*}
\chi_n\geq 0  , \  \int_{K}\chi_n(k)dk =1,\quad \textit{and}\quad \lim_{n\mapsto+\infty}\int_{K\backslash U}\chi_n(k)dk=0,
\end{eqnarray*}
for every neighbourhood $U$ of $e$ in $K$.\\
Let $F\in \mathcal{E}_{\lambda,p}^{*}(B(\mathbb{O}^2)$. For each $n$ define the function $F_n$ on $B(\mathbb{O}^2)$ by 
\begin{eqnarray*}
F_n(x)=\int_{K}\chi_n(k)F(k^{-1}.x)dk.
\end{eqnarray*}
Then  $ (F_n)_n$ converges pointwise to $F$ as $n$ goes to $+\infty$. Since $\Delta$ is $G-$invariant, then $F_n$ lies  in $\mathcal{E}_\lambda(B(\mathbb{O}^2))$.\\
For each $r\in[0,1[$ define a function $F^r$ on $K$, by $F^{r}(k)=F(rk.e_1)$.\\
Let $dk$ be the Haar measure of $K$. As 
\[\left(\int_{\partial B(\mathbb{O}^2)}\mid F(r\theta)\mid^2 d\theta\right)^{\frac{1}{2}}=\left(\int_{K}\mid F(rk.e_1)\mid^2dk\right)^{\frac{1}{2}},
\]
and noting that $F_n^r(k)=(\chi_{n}\ast F^r)(k)$, we get
\begin{eqnarray*}
\left(\int_{\partial B(\mathbb{O}^2)}\mid F_n(r\theta)\mid^2d\theta\right)^{\frac{1}{2}} & \le & \parallel \chi_n\parallel_{2}\parallel F^r\parallel_{1},\\
& \le & \parallel \chi_n\parallel_{2}\parallel F^r\parallel_{p},
\end{eqnarray*}
by the Haussedorf inequality. This shows that $(1-r^2)^{-\frac{\rho}{2}}\left(\int_{\partial B(\mathbb{O}^2)}\mid F_n(r\theta)\mid^2d\theta\right)^{\frac{1}{2}}$ is uniformly bounded for all $r\in [0,1[$. Thus, the  function $F_n$ lies in the space $\mathcal{E}_{\lambda,2}^{*}(B(\mathbb{O}^2))$. Therefore, for each $n$ there exists a function $f_n\in L^2(\partial B(\mathbb{O}^2))$ such that $ F_n= P_\lambda f_n$, by the if part of Theorem \ref{th1}. Moreover, we have
\begin{eqnarray*}
f_n(\omega)=\mid c(\lambda)\mid^{-2}\lim_{t\rightarrow +\infty}\frac{1}{t}\int_{B(0,t)} F_n(x)P_{-\lambda}(x,\omega)d\mu(x),
\end{eqnarray*}
in $L^{2}(\partial B(\mathbb{O}^2))$, by (\ref{invers}).\\
Next, denote by $g_n^t$ the function  
\begin{eqnarray*}
g_n^t(\omega)=\mid c(\lambda)\mid^{-2}\frac{1}{t}\int_{B(0,t)} F_n(x)P_{-\lambda}(x,\omega)d\mu(x).
\end{eqnarray*}
Below we will show that 
\begin{equation}\label{g_n}
\sup_{n\in \mathbb{N},t>0}\parallel g_n^{t}\parallel_{p} \le \mid c(\lambda)\mid^{-2}\gamma(\lambda,p)\parallel F\parallel_{*,p}.
\end{equation}
To do so, let $\phi\in L^{q}(\partial B(\mathbb{O}^2)$, with $\frac{1}{p}+\frac{1}{q}=1$. We have
\begin{eqnarray*}
 \left| \int_{\partial B(\mathbb{O}^2)}g_{n}^{t}(\omega)\overline{\phi(\omega)}d\omega\right|=
\end{eqnarray*}
\begin{eqnarray*}
 \frac{\mid c(\lambda)\mid^{-2}}{t}\left|\int_{\partial B(\mathbb{O}^2)}\left(\int_{B(0,t)} F_n(x)P_{-\lambda}(x,\omega)dx\right)\overline{\phi(\omega)}d\omega\right|.
\end{eqnarray*}
We may use the polar coordinates $x=r\theta$ ( $r\in[0,1[$ and $\theta\in \partial B(\mathbb{O}^2)$ ) to rewrite the right hand side of the above identity as
\begin{eqnarray*}
\frac{\mid c(\lambda)\mid^{-2}}{t} \left|\int_{\partial B(\mathbb{O}^2)}\left(\int^{\tanh(t)}_{0} \left(\int_{\partial B(\mathbb{O}^2)}P_{-\lambda}(r\theta,\omega)F_n(r\theta)d\theta \right) (1-r^2)^{-\rho-1}r^{15}dr\right)
\overline{\phi(\omega)}d\omega\right|.
\end{eqnarray*}
Since the Poisson kernel is symmetric in $\theta$ and $\omega$, then one can use the Fubini Theorem to show that 
\begin{eqnarray*}
 \left| \int_{\partial B(\mathbb{O}^2)}g_{n}^{t}(\omega)\overline{\phi(\omega)}d\omega\right|=
\end{eqnarray*}
\begin{eqnarray*}
 \mid c(\lambda)\mid^{-2}\frac{1}{t}\left|\int^{\tanh(t)} _{0}\left(\int_{\partial B(\mathbb{O}^2)}\overline{P_{\lambda}\phi}(r\theta)F_{n}(r\theta)d\theta \right)(1-r^2)^{-\rho-1}r^{15}dr\right|.
\end{eqnarray*}
By the H\"older inequality, the right side of the above equality is less than 
\begin{eqnarray*}
\frac{\mid c(\lambda)\mid^{-2}}{t} \int^{\tanh(t)} _{0}\left(
\int_{\partial B(\mathbb{O}^2)}\mid P_{\lambda}\phi(r\theta)\mid^{q} d\theta\right)^{\frac{1}{q}}\parallel F^r_{n}\parallel_{p}(1-r^2)^{-\rho-1}r^{15}dr.
\end{eqnarray*}
By the necessary condition already proved, we have
\begin{eqnarray*}
\left(
\int_{\partial B(\mathbb{O}^2)}\mid P_{\lambda}\phi(r\theta)\mid^{q} d\theta\right)^{\frac{1}{q}}\le \gamma(\lambda,q)(1-r^2)^{-\frac{\rho}{2}}\parallel \phi\parallel_{q},
\end{eqnarray*} 
and using
\begin{eqnarray*}
\parallel F^r_{n}\parallel_{p} \le \parallel F^{r}\parallel_{p}
\le (1-r^2)^{\frac{\rho}{2}}\parallel F\parallel_{\ast,p},
\end{eqnarray*} 
we deduce  that  
\begin{eqnarray*}
\left| \int_{\partial B(\mathbb{O}^2)}g_{n}^{t}(\omega)\overline{\phi(\omega)}d\omega\right|\le \mid c(\lambda)\mid^{-2} \gamma(\lambda,p)\parallel \phi\parallel_{q}\parallel F\parallel_{\ast,p},
\end{eqnarray*}
for every $t>0$ and $n$.\\
Next, taking the supermum over all $\phi$ such that $\parallel \phi\parallel_{q}=1$ we get the estimate (\ref{g_n}).\\
Now for $q$  a positive number with $\frac{1}{p}+\frac{1}{q}=1$, let  $T_n$ be the  linear form  on $L^{q}(\partial B(\mathbb{O}^2))$ defined by:
\[
T_n(\phi)=\int_{\partial B(\mathbb{O}^2)}f_n(\omega)\overline{\phi(\omega)}d\omega.
\]
Let $\phi$ be a $\mathbb{C}-$valued continuous function on $\partial B(\mathbb{O}^2)$ we have:
\begin{eqnarray*}
 T_n(\phi)=\lim_{t\mapsto+\infty} \int_{\partial B(\mathbb{O}^2)}g_n^t(\omega)\overline{\phi(\omega)}d\omega.
\end{eqnarray*}
Using on one hand  H\"older's inequality and on the other hand the estimate (\ref{g_n}), we see that 
\begin{eqnarray*}
\left| \int_{\partial B(\mathbb{O}^2)}g_n^t(\omega)\overline{\phi(\omega)}d\omega\right| \le \mid c(\lambda)\mid^{-2}\gamma(\lambda,p)\parallel F\parallel_{*,p}\parallel \phi \parallel_{q}.
\end{eqnarray*}
Next, taking the supremuum over all continuous function $\phi$ with $\parallel \phi \parallel_{q}=1$ in the above inequality, we deduce that the linear functionals $T_n$ are uniformly bounded, with 
\begin{equation}
\sup_n\parallel T_n\parallel\le \mid c(\lambda)\mid^{-2}\gamma(\lambda,p)\parallel F\parallel_{*,p},
\end{equation} 
where $\parallel .\parallel$ stands for the operator norm.\\
Thanks to the Banach-Alaoglu-Bourbaki Theorem, there exists a subsequence of bounded operators $(T_{n_j})_j$ which converges to a bounded linear operator $T$ on $L^{q}(\partial B(\mathbb{O}^2))$, under the $\ast$-weak topology, with $\parallel T\parallel\le \mid c(\lambda)\mid^{-2}\gamma(\lambda,p)\parallel F\parallel_{*,p}$.\\
By the Riesz representation theorem, there exists a unique function $f\in L^{p}(\partial B(\mathbb{O}^2))$ such that 
$T(\phi)=\int_{ \partial B(\mathbb{O}^2)}f(\omega)\overline{\phi(\omega)}d\omega$ with $\parallel f\parallel_{p}=\parallel T\parallel$. Therefore
\begin{equation}\label{ineq4}
\parallel f\parallel_{p}\le \mid c(\lambda)\mid^{-2}\gamma(\lambda,p)\parallel F\parallel_{*,p}.
\end{equation}
Now, let $\phi_{x}(\omega)=P_{\lambda}(x,\omega)$. Then, $T_{n}(\phi_{x})(\omega)=F_{n}(x)$. Since, on one hand  $\lim_{n\rightarrow +\infty}F_{n}(x)=F(x)$ and on the other hand  $\lim_{j\rightarrow +\infty} T_{n_j}(\phi_{x})=T(\phi_{x})$, we get $F=P_{\lambda}f$. The estimate $\parallel f\parallel_{p} \le \mid c(\lambda)\mid^{-2}\gamma(\lambda,p)\parallel P_{\lambda}f\parallel_{*,p}$ follows from (\ref{ineq4}). This finishes the proof of the main Theorem \ref{th2}.
\section{ Proof of the Key Lemma.}
This section is devoted to the proof of the Key Lemma of this paper. Recall that 
\begin{equation} \label{szego intg}
[\Psi_r(\lambda)f](\theta)= \int_{\partial B(\mathbb{O}^2)}\Psi_{r}(\lambda,\theta,\omega)f(\omega)d\omega,
\end{equation}
where $ \Psi_{r}(\lambda,\theta,\omega)=\Psi(r\theta,\omega)^{\frac{-i\lambda-\rho}{2}}$.\\
To prove the uniform $L^2-$uniform boundedness in $r\in[0,1[$ of the family of operators $(\Psi_r(\lambda))_{r\in[0,1[}$, we will adapt on $\partial B(\mathbb{O}^2)$, in a uniform manner in $r\in[0,1[$, the method of proving the $T(1)-$Theorem of David-Journ\'e and Semens \cite{D}, for this we follow the program accomplished by Y.Meyer in his new proof of the $T(1)-$Theorem in $L^{2}(\mathbb{R}^{n})$ which is based on the Cotlar-Stein Lemma.\\
To do this, we endow $\partial B(\mathbb{O}^2)$ with a non-isotropic metric $d$, so that $(\partial B(\mathbb{O}^2),d)$ becomes a space of homogeneous type in the sense of Coifman and Weiss, see \cite{C}.\\
More precisely, for $a\in \overline{B(\mathbb{O}^2)}, b\in \overline{B(\mathbb{O}^2)}$, we define 
\begin{eqnarray*}
d(a,b)=[1-2<a,b>_{\mathbb{R}}+\Phi(a,b)]^{\frac{1}{4}}.
\end{eqnarray*}
Since, $\Phi(a,b)$ is $O_{\mathbb{O}}(2)-$invariant and $O_{\mathbb{O}}(2)\subset O(16)$, we have 
\begin{eqnarray*}
d(h.a,h.b)=d(a,b),
\end{eqnarray*}
for every $h\in O_{\mathbb{O}}(2)$.\\
According to (\ref{Psi=1-[]}) the non-isotropic metric $d(a,b)$ may also be written as 
\begin{eqnarray*}
d(a,b)=\mid 1-[a,b]\mid^{\frac{1}{2}}.
\end{eqnarray*}
\begin{prop}\label{metric}
i) The triangle inequality
\begin{eqnarray*}
  d(a,c)\leq d(a,b)+d(b,c),
\end{eqnarray*}
holds for all $a,b\ and\ c \ in\  \overline{B(\mathbb{O}^2)}$.\\ 
ii) $d$ is a metric on $\partial B(\mathbb{O}^2)$.\\
iii) The volume of $B(\omega,\delta)$ with respect to the superficial measure of $\partial B(\mathbb{O}^2)$ behaves as $\delta^{2\rho}$, where $\rho=11$ and  $B(\omega,\delta)=\{\theta \in \partial B(\mathbb{O}^2) :d(\theta,\omega)<\delta\}$.
\end{prop}
\textbf{Proof.}
i) Put $e_2=(0,1)$. Since $d$ is bi-invariant by $O_{\mathbb{O}}(2)$, we may take $b=re_2$ ( $0\le r \le 1$ ) and we then have to prove 
\begin{equation}\label{ineq1}
d(a,c)\le \mid 1-ra_2\mid^{\frac{1}{2}}+\mid 1-rc_2\mid^{\frac{1}{2}}.
 \end{equation}
Notice that (\ref{ineq1}) is obvious if $c_2=0$.\\
If $c_2\neq 0$, then one has to prove that 
\begin{eqnarray*}
|1-(\overline{a_1}c_2)(c_2^{-1}c_1)-a_2\overline{c_2}|\le \left[\mid 1-ra_2\mid^{\frac{1}{2}}+\mid 1-rc_2\mid^{\frac{1}{2}}\right]^{2}.  
\end{eqnarray*}
The left hand side of the above inequality is less than 
\begin{eqnarray*}
|1-a_2\overline{c_2}|+\mid a_1\mid\mid c_1\mid.
\end{eqnarray*}
Next, since
\begin{eqnarray*}
|1-a_2\overline{c_2}|\le |1-ra_2|+|1-rc_2|,
\end{eqnarray*}
and 
\begin{eqnarray*}
\mid a_1\mid^{2}\le 1-\mid a_2\mid^2\le 1-r^2\mid a_2\mid^2 \le 2\mid 1-ra_2\mid,
\end{eqnarray*}
and similar estimate for $\mid c_1\mid$, then (\ref{ineq1}) holds.\\
ii) To prove that $d$ is a metric on $\partial B(\mathbb{O}^2)$, we have only to show that $d(a,b)=0$ if and only if $a=b$.\\
By the $O_{\mathbb{O}}(2)-$invariance of $d$ it suffices to prove it for $b=e_2$, which is obvious.\\
For the proof of iii) we will need the following standard calculus Lemma: 
\begin{lem}\label{int}
Let $f$ be a $\mathbb{C}-$valued function on $\partial B(\mathbb{O}^2)$ such that \\$f(\omega_1,\omega_2)=g(\omega_1)$.
Then, we have
\begin{eqnarray*}
\int_{\partial B(\mathbb{O}^2)}f(\omega)d\omega=c  \int_{\{x\in\mathbb{O}:|x|<1\}}g(x)(1-|x|^{2})^3 dm(x).
\end{eqnarray*}
\end{lem}
Since the metric $d$ is $O_{\mathbb{O}}(2)-$invariant  as well as the superficial measure $d\theta$, we have $V(B(\omega,\delta))=V(B(e_1,\delta))$. Therefore
\begin{eqnarray*}
  V(B(\omega,\delta))= \int_{ \{ \theta\in\mathbb{O}^2, |1-\theta_1|<\delta^2  \}}d\theta.
\end{eqnarray*}
Next, use Lemma \ref{int} to get:
\[
  V(B(\omega,\delta))=c  \int_{\{x\in \mathbb{O}: |x|<1; |1-x|< \delta^2\}}(1-|x|^2)^3 dm(x).
\]
Put $1-x=t(\cos(\alpha)+\sin(\alpha)y)$, where $t>0$, $\alpha \in [0,\pi]$ and $y\in \mathbb{O}$ such that $\mathfrak{Re}(y)=0$ and $|y|=1$. Then the above integral may be rewritten as:
  \begin{eqnarray*}
  V(B(\omega,\delta))
  = c  \int_{\{(\alpha,t) \in [0,\pi]\times ]0,\delta^2[;\mid 1-te^{i\alpha})\mid<1 \}}(2\cos(\alpha)-t)^3 t^{10} \sin^{6}(\alpha) dt d\alpha,
\end{eqnarray*}
from which we deduce easily that $V(B(\omega,\delta))\le c  \delta^{22}$. This finishes the proof of Proposition \ref{metric}.\\ 
Next, to make  Meyer program works in our case, we follow the same line as in the proof in \cite{BI}.\\
\textbf{Step 1:Uniform Calderon-Zygmund type estimates}
\begin{prop}\label{calderon}
There exists a positive constant $c $ such that the following estimates hold\\

i)\begin{eqnarray*}
 \quad  \sup_{0\leq r<1}|\Psi_r(\lambda,\theta,\omega)|\leq c  d(\theta,\omega)^{-2\rho},
 \end{eqnarray*}
for every $\theta$ and $\omega$ in $\partial B(\mathbb{O}^2)$.\\
ii)
\begin{eqnarray*}
 \sup_{0\leq r<1}|\Psi_r(\lambda,\theta,\omega)-\Psi_r(\lambda,\theta^\prime,\omega)|\leq c (1+\mid\lambda\mid) \frac{d(\theta,\theta^\prime)}{d(\theta,\omega)^{2\rho+1}}, 
\end{eqnarray*}
for every $\theta, \theta^\prime, \omega \ in\  \partial B(\mathbb{O}^2)$ such that $d(\theta,\omega)\geq 2d(\theta, \theta^\prime)$.\\
iii) 
\begin{eqnarray*}
\sup_{0\le r<1}\left| \int_{d(\theta,\omega)\le \delta}\Psi_{r}(\lambda,\theta,\omega)d\omega\right| \le c (1+\frac{1}{\mid \lambda\mid}),
\end{eqnarray*}
for every $\delta > 0$.
\end{prop}
\textbf{Proof.} Notice that 
\begin{eqnarray*}
\Psi_r(\lambda,\theta,\omega)=|1-r[\theta,\omega]|^{-i\lambda-\rho}.
\end{eqnarray*}
i) We have 
\begin{eqnarray*}
|1-[\theta,\omega]|=\mid 1-r[\theta,\omega]-(1-r)[\theta,\omega]\mid \le \mid 1-r[\theta,\omega]\mid +(1-r)|[\theta,\omega]|.
\end{eqnarray*}
Since $|[\theta,\omega]|\leq 1$, 
and for all $r\in [0,1[$, $ 1-r\le\mid 1-r[\theta,\omega]\mid$ , then we get
\begin{equation}\label{ineq5}
 \mid 1-r[\theta,\omega]\mid^{-1}\le 2 \mid 1-[\theta,\omega]|^{-1},\forall r\in [0,1[
\end{equation}
and (i) follows.\\
ii) Now, by the mean calculus Lemma  we obtain
\begin{eqnarray*}
\mid\Psi_r(\lambda,\theta,\omega)-\Psi_r(\lambda,\theta^{\prime},\omega)\mid\le \frac{|i\lambda+\rho|}{\rho}\left||1-r[\theta,\omega]|^{-\rho}-|1-r[\theta^\prime,\omega]|^{-\rho}\right|. 
\end{eqnarray*} 
Then the  proof of (ii) will be based on the identity
\begin{eqnarray*}
|1-r[\theta,\omega]|^{-\rho}-|1-r[\theta^\prime,\omega]|^{-\rho}=\sum\limits_{j=0}^{\rho-1}\frac{\mid 1-r[\theta^{\prime},\omega]\mid - \mid 1-r[\theta,\omega]\mid}{\mid 1-r[\theta,\omega]\mid^{j+1}\mid 1-r[\theta^{\prime},\omega]\mid^{\rho-j}}.
\end{eqnarray*}
By (\ref{ineq5}), we have
\begin{eqnarray*}
\mid 1-r[\theta,\omega]\mid^{-1}\le 2 \mid 1-[\theta,\omega]|^{-1},
\end{eqnarray*} 
and 
\begin{eqnarray*}
\mid 1-r[\theta^{\prime},\omega]\mid^{-1}\le 2 \mid 1-[\theta^{\prime},\omega]\mid^{-1}\le 8\mid 1-[\theta,\omega]\mid^{-1}.
\end{eqnarray*}
The last inequality is a consequence of the triangle inequality and the hypotheses $d(\theta,\omega)\geq 2d(\theta, \theta^\prime)$. It is clear from above that 
\begin{eqnarray*}
\mid\Psi_r(\lambda,\theta,\omega)-\Psi_r(\lambda,\theta^{\prime},\omega)\mid\le c  \frac{(\lambda^2+\rho^2)^{\frac{1}{2}}}{\rho} \frac{\mid[\theta-\theta^{\prime},\omega]\mid}{\mid 1-[\theta,\omega]\mid^{\rho+1}},
\end{eqnarray*}
for some numerical constant $c $.\\
We claim that 
\begin{equation}\label{ineq2}
|[\theta-\theta^\prime, \omega]|\leq d(\theta,\theta^\prime)(d(\theta,\theta^\prime)+2d(\theta,\omega)),
\end{equation}
for any $\theta, \theta^\prime \ and \ \omega$ in $\partial B(\mathbb{O}^2)$.\\
In the proof of (\ref{ineq2})  we may take   $\theta=e_1=(1,0)$ ( by the $O_{\mathbb{O}}(2)-$invariance of the non-isotropic distance ) and one have to prove
\begin{equation}\label{ineq3}
|[e_1-\theta^\prime,\omega]|\le \mid 1-\theta_1^{\prime}\mid^{\frac{1}{2}}(\mid 1-\theta_1^{\prime}\mid^{\frac{1}{2}}+2\mid 1-\omega_1\mid^{\frac{1}{2}}).
\end{equation}
It is obvious for $\omega_2=0$.\\
If $\omega_2\neq 0$, then 
\begin{eqnarray*}
\mid [e-\theta^\prime,\omega]\mid
 &=& \mid((1-\overline{\theta^\prime_1})\omega_2)(\omega_2^{-1}\omega_1)-\theta^\prime_2\overline{\omega_2}\mid\\
&\leq& |1-\overline{\theta^\prime_1}|+|\theta^\prime_2||\omega_2|.
\end{eqnarray*}
Next from $|\theta^\prime_2|^2=1-|\theta^\prime_1|^2\leq2|1-\theta^\prime_1|$ and  similar estimate for $\mid \omega_2\mid$, we get (\ref{ineq3}). This finishes the proof of (\ref{ineq2}).\\
Collecting the above results we can conclude that 
\begin{eqnarray*}
\mid\Psi_r(\lambda,\theta,\omega)-\Psi_r(\lambda,\theta^{\prime},\omega)\mid\le c  (1+|\lambda|) \frac{d(\theta,\theta^{\prime})}{d(\theta,\omega)^{2\rho+1}}.
\end{eqnarray*}
This finishes the proof of (ii).\\
iii) By the $O_{\mathbb{O}}(2)-$invariance of the metric $d$ and the measure $d\omega$, we have :
\begin{equation}\label{K-invar}
\int_{d(\omega,\theta)<\delta}\Psi_r(\lambda,\theta,\omega)d\omega=\int_{|1-\omega_1|<\delta^2} |1-r\omega_1|^{-i\lambda-\rho}d\omega.
\end{equation}
It is clear that (\ref{K-invar}) is uniformly bounded for every $r\in[0,\dfrac{1}{2}[$ and $\delta>0$. 
\\
To show the uniform boundedness for  $r\in[\dfrac{1}{2},1[$ we use Lemma \ref{int} to get 
\[
\int_{d(\omega,\theta)<\delta}\Psi_r(\lambda,\theta,\omega)d\omega=\int_{\lbrace x\in \mathbb{O}: |x|<1, |1-x|<\delta^2\rbrace}|1-rx|^{-i\lambda-\rho}\left( 1-|x|^2\right)^3 dm(x).
\]
We put $1-rx=t(\cos(\alpha)+\sin(\alpha)y)$ where $t>0$, $\alpha \in [0,\pi]$ and $y\in \mathbb{O}$ such that $\mathfrak{Re}(y)=0$ and $|y|=1$. The above integral may be rewritten 
\[
\int_{d(\omega,\theta)<\delta}\Psi_r(\lambda,\theta,\omega)d\omega=r^{-14}\int_{\Gamma_{r,\delta}} t^{-i\lambda-4}\mid r^2 -|te^{i\alpha}-1|^2\mid^3 \sin^6(\alpha)dt d\alpha,
\]
where $\Gamma_{r,\delta}$ is the set of $(t,\alpha) \in ]0,\infty[ \times [0,\pi]$ such that $|te^{i\alpha}-1|<r$, and $|(te^{i\alpha}-(1-r)|<r\delta^2$.\\
Next, replacing $\sin^6(\alpha)$ by
\[
\sin^6(\alpha)=\frac{1}{64}\left( 20-15(e^{2i\alpha}+e^{-2i\alpha})+6(e^{4i\alpha}+e^{-4i\alpha})-(e^{6i\alpha}+e^{-6i\alpha})\right) ,
\]
we are reduced to show the uniform boundedness of integrals of the following type
\[
\int_{\Gamma_{r,\delta}} t^{-i\lambda-4}\mid r^2 -|te^{i\alpha}-1|^2\mid^3 e^{ik\alpha}dt d\alpha.
\]
We may use our result in \cite{BI1}, to show that that there exists a positive constant $c $ such that :
\begin{eqnarray*}
\left| \int_{d(\omega,\theta)<\delta}\Psi_r(\lambda,\theta,\omega)d\omega \right| \leq c  (1+\dfrac{1}{|\lambda|}),
\end{eqnarray*}
for every $r\in [\frac{1}{2},1[$ and $\delta>0$. This finishes the proof of the proposition.\\
\textbf{Step 2: Uniform action of Szeg\"o-integrals on $\delta$-molecules.}\\
For $\eta>0$ and $0<\delta\leq 1$, we define the weight function 
\begin{eqnarray*}
\Omega_{\eta,\delta}(\theta,\omega)=\eta^{\delta}[\eta+d(\theta,\omega)]^{-\delta-2\rho},
\end{eqnarray*}
for every $\theta$ and $\omega$ in $\partial B(\mathbb{O}^2)$.\\
\textbf{Definition.}
 A $\mathbb{C}-$valued function $m$ on $\partial B(\mathbb{O}^2)$ is said to be a $\delta-$molecule centred at $\theta_{0}$ with width $\eta >0$ if $m$ satisfies
\begin{itemize}
\item[i)] $\mid m(\theta)\mid \leq \Omega_{\eta,\delta}(\theta,\theta_{0})$,
\item[ii)] $\mid m(\theta)-m(\theta')\mid \leq (\frac{d(\theta,\theta')}{\eta})^{\delta} \left( \Omega_{\eta,\delta}(\theta,\theta_{0}) + \Omega_{\eta,\delta}(\theta',\theta_{0})\right)$,
\item[iii)] $\int_{\partial B(\mathbb{O}^2)}m(\theta)d\theta=0$.
\end{itemize}
Let $M(\delta,\theta_0,\eta)$ denote the convex set of all $\delta-$molecules centered at $\theta_0$ with width $\eta$.\\
Then with the help of the above definition and using the uniform Calderon-Zygumund estimates on the Schwartz kernel of $\Psi_{r}(\lambda)$, we can prove 
\begin{prop}
The operator $\Psi_{r}(\lambda)$ transforms uniformly in $r\in[0,1[$, $\delta-$molecules into $\delta'-$molecules with $0<\delta'<\delta \leq 1$. More precisely, there exists a positive constant $c $ such that for any $m\in M(\delta,\theta_0,\eta)$ and any $0<\delta^{\prime}<\delta\le 1$, the function
\begin{eqnarray*}
m^{\prime}(\theta)=c (1+\mid \lambda\mid+\frac{1}{\mid \lambda\mid})(\Psi_{r}(\lambda)m)(\theta),  
\end{eqnarray*}
lies in $M(\delta^{\prime},\theta_0,\eta)$.
\end{prop}
\textbf{Step 3: Molecular resolution of $L^{2}(\partial B(\mathbb{O}^2))$}.\\
We built an adapted $\delta-$molecular resolution by means of the Poisson kernel $P(x,\omega)$. Namely, for $j=0,1,...$, we set 
\begin{eqnarray*}
\Delta_{j}(\theta,\omega)=P(\eta_{j+1}\theta,\omega)-P(\eta_{j}\theta,\omega),
\end{eqnarray*}
where $\eta_{j}=\frac{2(1-2^{-j})}{2^{-2j}+2(1-2^{-j})}$. 
\begin{prop}
Let $\delta>0$. Then we have\\
(i) The functions $\theta\rightarrow \Delta_{j}(\theta,\theta_{0})$ are $\delta-$molecules centered at $\theta_{0}$ with width $2^{-j}$.\\
(ii) $I_{H}=\sum\limits_{j\geq 0}\Delta_{j}$, where $H$ is the orthogonal of $1$ in $L^{2}(\partial B(\mathbb{O}^2))$.
\end{prop}
Now to get the estimate (\ref{Q1}) in the Key Lemma, we may  write the operator $\Psi_{r}(\lambda)$ as 
\begin{eqnarray*}
\Psi_{r}(\lambda)=\sum\limits_{j=0}^{+\infty}\Psi_{r,j}(\lambda),
\end{eqnarray*}
with $\Psi_{r,j}(\lambda)=\Psi_{r}(\lambda)\circ \Delta_{j}$.\\
Next, using the uniform action of $\Psi_{r}(\lambda)$ on molecules as well as the bounded mean property iii) in Proposition \ref{calderon}, we may apply uniformly in $ r\in [0,1[$ the Cotlar-Stein Lemma to obtain $||\Psi_{r}(\lambda)||_2\le c (1+\mid \lambda\mid+\frac{1}{\mid \lambda\mid}) $.\\
Finally, combining the estimate ii) of Proposition \ref{calderon}  as well as iii) of Proposition \ref{metric},  we get  the H\"ormander condition (\ref{Q2}) on the Schwartz kernel and the proof of the Key Lemma of this paper is finished.

\vspace{0.5cm}
\textsc{DEPARTMENT OF MATHEMATICS, FACULTY OF SCIENCES,\\
UNIVERSITY IBN TOFAIL, MOROCCO.}\\
E-mail address: a.boussejra@gmail.com, boussejra.abdelhamid@uit.ac.ma\\nadia.ourchane20@gmail.com

\end{document}